\documentclass[twoside,11pt]{article}
\usepackage{amsmath,mathtools,amsfonts,amscd,amsthm,graphics,latexsym,stmaryrd}
\usepackage[all]{xy}
\usepackage{mathptmx}
\usepackage{anysize}
\usepackage[bookmarks=true]{hyperref}
\usepackage{fancyhdr}

\marginsize{3.385cm}{3.385cm}{3cm}{3cm}

\newtheorem{thm}{Theorem}[section]

\newtheorem{Thm}[thm]{Theorem}
\newtheorem{cor}[thm]{Corollary}
\newtheorem{prop}[thm]{Proposition}
\newtheorem{lem}[thm]{Lemma}
\newtheorem{conj}[thm]{Conjecture}
    
\numberwithin{equation}{section}
\DeclareSymbolFont{oldsymbols}{OMS}{cmsy}{m}{n}
\DeclareSymbolFontAlphabet{\mathcal}{oldsymbols}

\let\jmath=\undefined
\DeclareSymbolFont{cmletters}{OML}{cmm}{m}{it}
\DeclareMathSymbol{\jmath}{\mathord}{cmletters}{"7C}

\let\coprod=\undefined
\DeclareSymbolFont{cmsymbols}{OMS}{cmsy}{m}{n}
\DeclareSymbolFont{cmlargesymbols}{OMX}{cmex}{m}{n}
\DeclareMathSymbol{\coprod}{\mathop}{cmlargesymbols}{"60}

\makeatletter 
\def\mathcenterto#1#2{\mathclap{\phantom{#1}\mathclap{#2}}\phantom{#1}} 
\let\old@wt\widetilde 
\def\widetildeto#1#2{\mathcenterto{#2}{\old@wt{\mathcenterto{#1}{#2}}}} 
\makeatother 
 
\def\wt{\widetildeto{o}} 

\newcommand{\Hilb}{\mathrm{Hilb}}
\newcommand{\Sym}{\mathrm{Sym}}
\newcommand{\Ast}{A^*}
\newcommand{\AT}{A_{\bT}^*}

\newcommand{\Aorb}{A_{\mathrm{orb}}^*}
\newcommand{\ATorb}{A_{\bT, \mathrm{orb}}^*}

\newcommand{\wtd}{\widetilde}
\newcommand{\qcr}{\cup_{\mathrm{qc}}}

\newcommand{\vir}{\mathrm{vir}}
\newcommand{\mov}{\mathrm{mov}}
\newcommand{\Hom}{\mathrm{Hom}}

\newcommand{\st}{\mathrm{st}}
\newcommand{\M}{\overline{M}}
\newcommand{\N}{\overline{M}_{0,3}}
\newcommand{\I}{\overline{I}}
\newcommand{\age}{\mathrm{age}}
\newcommand{\Aut}{\mathrm{Aut}}
\newcommand{\eT}{e_{\mathbb{T}}}
\newcommand{\ev}{\mathrm{ev}}
\newcommand{\Sn}{S^{[n]}}
\newcommand{\spa}{\textrm{ }}

\newcommand{\Ar}{\mathcal{A}_{r}}
\newcommand{\An}{\mathcal{A}_{n}}
\newcommand{\bae}{\bar{e}}

\newcommand{\bC}{\mathbb{C}}

\newcommand{\bP}{\mathbb{P}}
\newcommand{\bQ}{\mathbb{Q}}
\newcommand{\bT}{\mathbb{T}}
\newcommand{\bZ}{\mathbb{Z}}

\newcommand{\cC}{\mathcal{C}}
\newcommand{\cI}{\mathcal{I}}
\newcommand{\cO}{\mathcal{O}}
\newcommand{\cP}{\mathcal{P}}
\newcommand{\cQ}{\mathcal{Q}}
\newcommand{\cV}{\mathcal{V}}
\newcommand{\cX}{\mathcal{X}}

\newcommand{\fa}{\mathfrak{a}}
\newcommand{\fB}{\mathfrak{B}}
\newcommand{\fM}{\mathfrak{M}}
\newcommand{\fm}{\mathfrak{m}}
\newcommand{\fp}{\mathfrak{p}}
\newcommand{\fS}{\mathfrak{S}}
\newcommand{\fz}{\mathfrak{z}}

\newcommand{\al}{\alpha}
\newcommand{\be}{\beta}
\newcommand{\g}{\gamma}
\newcommand{\G}{\Gamma}
\newcommand{\de}{\delta}
\newcommand{\si}{\sigma}
\newcommand{\Si}{\Sigma}
\newcommand{\la}{\lambda}
\newcommand{\La}{\Lambda}
\newcommand{\io}{\iota}
\newcommand{\jm}{\jmath}
\newcommand{\vect}{\overrightarrow}

\def\co{\colon\thinspace}

\begin{document}

\title{\Large \bf Strengthening the cohomological crepant resolution conjecture for Hilbert--Chow morphisms}

\author{Wan Keng Cheong}
\date{}

\maketitle

\begin{abstract}
Given any smooth toric surface $S$, we prove a SYM-HILB correspondence which relates the 3-point, degree zero, extended Gromov--Witten invariants of the $n$-fold symmetric product stack $[\Sym^n(S)]$ of $S$ to the 3-point extremal Gromov--Witten invariants of the Hilbert scheme $\Hilb^n(S)$ of $n$ points on $S$. As we do not specialize the values of the quantum parameters involved, this result proves a strengthening of Ruan's Cohomological Crepant Resolution Conjecture for the Hilbert--Chow morphism $\Hilb^n(S) \to \Sym^n(S)$ and yields a method of reconstructing the cup product for $\Hilb^n(S)$ from the orbifold invariants of $[\Sym^n(S)]$.

\end{abstract}

\setcounter{section}{-1}

\section{Introduction}
\subsection{Overview}
Let $S$ be a smooth complex surface and $n$ a positive integer. The symmetric group $\fS_n$ on $n$ letters acts on the $n$-fold product $S^n$ by
$$ g\cdot(s_1,\ldots,s_n)=(s_{g(1)},\ldots,s_{g(n)}).$$ 
The quotient scheme $S^n/\fS_n$, denoted by $\Sym^n(S)$, is referred to as the $n$-fold symmetric product of $S$. It is singular for $n \geq 2$, but the quotient stack $[S^n/\fS_n]$ (cf. Sect. \ref{s: DefSym}) is a (smooth) orbifold for every $n$. We denote $[S^n/\fS_n]$ by $[\Sym^n(S)]$ and call it the $n$-fold symmetric product stack of $S$. Note that $\Sym^n(S)$ is the coarse moduli scheme of $[\Sym^n(S)]$.

The Hilbert scheme of $n$ points on $S$, written as $\Hilb^n(S)$ or $S^{[n]}$, parametrizes zero-dimensional closed subscheme $Z$ of $S$ satisfying $\mathrm{dim}_{\bC} \, H^0(Z,\cO_Z)=n$. 
Moreover, there exists a resolution of singularities
$\rho \co \Hilb^n(S)\to \Sym^n(S)$
defined by
$$\rho(Z)=\sum_{p\in S}\ell(\cO_{Z,p}) \, [p],$$
where $\ell(\cO_{Z,p})$, the length of $\cO_{Z,p}$, is simply the multiplicity of $p$ in $Z$; see Fogarty \cite{Fo}. The resolution $\rho$ is called the Hilbert--Chow morphism and is also crepant (see Beauville \cite{Bea}), i.e., 
$$K_{\Hilb^n(S)}=\rho^* K_{\Sym^n(S)},$$
where $K_{\Hilb^n(S)}$ (resp. $K_{\Sym^n(S)}$) is the canonical class of $\Hilb^n(S)$ (resp. $\Sym^n(S)$).
Furthermore, Fu and Namikawa \cite{FuN} show that $\rho$ provides a unique crepant resolution for $\Sym^n(S)$.

The following diagram
$$
\xymatrix{&[\Sym^n(S)]  \ar[d]^{c}  \\ \Hilb^n(S) \ar[r]^{\rho} &\Sym^n(S)}
$$
summarizes the relationships among these spaces. Here $c$ is the canonical map to the coarse moduli space.

Theoretical physicists believe that {\it string theory on an orbifold and string theory on any crepant resolution should belong to the same family}. As for the above examples, it is expected that there is an equivalence between string theories of $[\Sym^n(S)]$ and $\Hilb^n(S)$. 

The physical principle has led to various mathematical predictions; see, for example, Ruan \cite{R}, Bryan--Graber \cite{BrG}, Coates--Iritani--Tseng \cite{CoIT}, and Coates--Ruan \cite{CoR}. In this article, we are particularly interested in the following conjecture, which is referred to as the Cohomological Crepant Resolution Conjecture (abbreviated as CCRC).

\begin{conj}[\cite{R}] \label{R-CCRC}
Let $\cX$ be a (smooth) Gorenstein orbifold and $X$ its coarse moduli space. Suppose that $X$ admits a crepant resolution $Y$. Then the Chen--Ruan cohomology ring of $\cX$ is isomorphic to the quantum corrected cohomology ring of $Y$.
\end{conj}

The notions of Chen--Ruan cohomology and quantum corrected cohomology will be recalled later.

There are several examples for which CCRC is known to be true. For instance, Fantechi and G\"ottsche \cite{FaG}, and independently Uribe \cite{U} apply the results of Lehn and Sorger \cite{LeS} (see also Lehn's work \cite{Le} and Z. Qin and W. Wang's work \cite{LQW1}) to establish the validity of CCRC for $\Sym^n(S)$, where $S$ is a smooth complex projective surface with trivial canonical class. Note that a different proof is also obtained by Qin and Wang \cite{QW}. Moreover, when $S$ is a smooth, simply-connected, complex projective surface, J. Li and W.-P. Li \cite{LL} determine the 2-point extremal Gromov--Witten invariants of $\Hilb^n(S)$ and confirm that there is a linear isomorphism between the Chen--Ruan cohomology ring of $[\Sym^n(S)]$ and the quantum corrected cohomology ring of $\Hilb^n(S)$ which respects multiplication by divisor classes\footnote{Recently, W.-P. Li and Z. Qin \cite{LQ} have proved that CCRC is true in that case, i.e., the linear isomorphism is indeed an algebra isomorphism. Their proof relies on the validity of CCRC for symmetric products of smooth toric surfaces, which is confirmed in this paper; see Corollary \ref{CCRC}.}. However, CCRC in the case of $\Sym^n(S)$, for an arbitrary smooth toric surface $S$, has not yet been fully verified. This case will be the main focus of this paper.

Let $\bT=(\bC^{\times})^2$. The $\bT$-equivariant cohomology of a point is simply the polynomial algebra in two variables $t_1, t_2$. In this article, we assume that all equivariant cohomology groups are with respect to the torus $\bT$.

Let $S$ be a smooth toric surface. The basic objects of this paper are 3-point, degree zero, extended Gromov--Witten invariants of $[\Sym^n(S)]$ and 3-point extremal Gromov--Witten invariants of $\Hilb^n(S)$. They are encoded in what we call the extended 3-point functions
$$\langle -, -, - \rangle^{[\Sym^n(S)]}(u) \in \bQ(t_1, t_2)[[u]] \spa \spa \textrm{ (cf. Sect. \ref{s: OrbGW})}$$ 
and the extremal 3-point functions
$$\langle -, -, - \rangle^{\Hilb^n(S)}(q) \in \bQ(t_1, t_2)[[q]] \spa \spa \textrm{ (cf. Sect. \ref{s: extrGW})}$$
respectively.

Our goal is to construct a SYM-HILB correspondence, and to prove a strengthening of CCRC for the Hilbert--Chow morphism $\rho \co \Hilb^n(S) \to \Sym^n(S)$. The following is our main result.

\begin{Thm}\label{main}
Let $q=-e^{iu}$ where $i$ is a square root of $-1$. For any smooth toric surface $S$ and any positive integer $n$, there is an isometric isomorphism $L$ which maps the equivariant Chen--Ruan cohomology of $[\Sym^n(S)]$ onto the equivariant cohomology of $\Hilb^n(S)$, and which satisfies the identities
$$
\langle \al_1, \al_2, \al_3 \rangle^{[\Sym^n(S)]}(u)=\langle L(\al_1), L(\al_2), L(\al_3) \rangle^{\Hilb^n(S)}(q)
$$
for any equivariant Chen--Ruan cohomology classes $\al_1, \al_2, \al_3$. 
\end{Thm}

We will make the correspondence $L$ explicit in Sect. \ref{s:3}. Roughly speaking, it maps the fixed-point basis for the equivariant Chen--Ruan cohomology of $[\Sym^n(S)]$ to the Nakajima basis for the equivariant cohomology of $\Hilb^n(S)$. Our proof uses a localization technique and relies on the case of $S=\bC^2$, which follows from the results of Bryan--Graber \cite{BrG} and Okounkov--Pandharipande \cite{OP}.

The equivariant Chen--Ruan cohomology ring of the symmetric product stack $[\Sym^n(S)]$ is given by extended 3-point functions $\langle -, -, - \rangle^{[\Sym^n(S)]}(u)$ with $u$ being specialized to $0$, while the equivariant quantum corrected cohomology ring of the Hilbert scheme $\Hilb^n(S)$ is defined by extremal 3-point functions $\langle -, -, - \rangle^{\Hilb^n(S)}(q)$ with $q$ being set to $-1$. Thus, it is an immediate consequence of Theorem \ref{main} that CCRC is valid for $[\Sym^n(S)]$ and $\Hilb^n(S)$.
\begin{cor}\label{CCRC}
The equivariant Chen--Ruan cohomology ring of $[\Sym^n(S)]$ is isomorphic to the equivariant quantum corrected cohomology ring of $\Hilb^n(S)$. 
\end{cor}

On the other hand, by taking $q =0$, we have the following.
\begin{cor}\label{cup}
The cup product for $\Hilb^n(S)$ can be recovered from the extended 3-point functions of $[\Sym^n(S)]$.
\end{cor}

If we have closed-form formulas for the symmetric product orbifold invariants involved, the cup product for the Hilbert scheme of points can be written down explicitly.  

Furthermore, we can even use the map $L$ and the setting of this paper to compare the full Gromov--Witten theories of $[\Sym^n(\Ar)]$ and $\Hilb^n(\Ar)$, where $\Ar$ is the minimal resolution of the quotient variety $\bC^2/\mu_{r+1}$ ($\mu_{r+1}$ is the group of $(r+1)$-th roots of unity). We refer the reader to Cheong--Gholampour \cite{ChG} and Maulik--Oblomkov \cite{MO} for more details.

\subsection{Outline}
We investigate the extended 3-point functions of $[\Sym^n(S)]$ in Sect. \ref{s: Sym} and the extremal 3-point functions of $\Hilb^n(S)$ in Sect. \ref{s: Hilb}. In Sect. \ref{s:3}, we give a concrete description of the SYM-HILB correspondence mentioned above and use the results of Sect. \ref{s: Sym} and Sect. \ref{s: Hilb} to show Theorem \ref{main}.

\subsection{Setting}
Throughout the paper, we let $S$ be a smooth toric surface acted upon by the torus $\bT=(\bC^{\times})^2$.

The surface $S$ is determined by a fan $\Si$ which is a finite collection of strongly convex rational polyhedral cones $\si$ contained in $N=\Hom(M, \bZ)$, where $M \cong \bZ^2$. That is, $S$ is obtained by gluing together affine toric varieties $S_{\si}$ and $S_{\tau}$ along $S_{\si\cap\tau}$ for $\si, \tau \in \Si$. Here, for example, the coordinate ring of $S_\si$  is $\bC[{\si}^{\vee}\cap M]$, which is the $\bC$-algebra with generators $\chi^m$ for $m \in {\si}^{\vee}\cap M$ and multiplication defined by $\chi^m \chi^{m^\prime}=\chi^{m+m^\prime}$. Note that ${\si}^{\vee}\cap M$ is, by definition, the set of elements $m \in M$ satisfying $v(m) \geq 0$ for all $v \in \si$. It is a finitely generated semigroup, and so the $\bC$-algebra $\bC[{\si}^{\vee}\cap M]$ is finitely generated.

In addition, $S$ has finitely many $\bT$-invariant subvarieties, and so it has a finite number of $\bT$-fixed points, denoted by 
$$x_1,\ldots,x_s.$$ 
(We do not study smooth toric surfaces without $\bT$-fixed points as they are not interesting in equivariant theory.)
 
For each $i$, the point $x_i$ lies on
$$U_i:=S_{\si_i}$$ 
for some $\si_i\in \Si$. As $S$ is smooth and $U_i$ possesses a unique $\bT$-fixed point $x_i$, we see that $U_i$ must be isomorphic to the affine plane with $x_i$ corresponding to the origin (cf. Fulton \cite{Ful}). However, $S$ is not necessarily the union $\bigcup_{i=1}^s U_i$.

From here on, let us fix the above setting on the open sets $U_i$ and the fixed points $x_i$. We denote by
$$L_i \hspace{3mm} \textrm{and} \hspace{3mm} R_i$$ 
the tangent weights at $x_i$.

For ease of exposition, we need some other notation. Below is the one that will be used frequently.

\paragraph{\underline{\bf Notation:}}
\begin{enumerate}
\item To avoid doubling indices, we identify 
$$A^i(X)=H^{2i}(X; \bQ), \hspace{3mm} A_i(X)=H_{2i}(X; \bQ), \hspace{3mm}\textrm{and } \hspace{3mm} A_i(X; \bZ)=H_{2i}(X; \bZ),$$ 
just to name a few, for any complex variety $X$ to appear in this article (note that we drop $\bQ$ but not $\bZ$). They will be referred to as cohomology or homology groups rather than Chow groups.

\item
\begin{enumerate}
\item Denote by $t_1, t_2$ the generators of the equivariant cohomology of a point, i.e., $\AT(\textrm{point})=\bQ[t_1, t_2]$.
\item $V_{\fm}=V\otimes_{\bQ[t_1,t_2]}\bQ(t_1,t_2)$ for each $\bQ[t_1,t_2]$-module $V$.
\end{enumerate}

\item For any space $X$ with a $\bT$-action, $X^\bT$ denotes the $\bT$-fixed locus of $X$.

\item An orbifold $\cX$ is a smooth Deligne--Mumford stack of finite type over $\bC$. Denote by $$c\co \cX \to X$$ the canonical map to the coarse moduli space.

\item For any object $\cO$, $\cO^n$ means that $\cO$ repeats itself $n$ times.

\item Let $\si$ be a partition of a nonnegative integer.
 \begin{enumerate}
 \item $\ell(\si)$ denotes the length of $\si$. Unless otherwise stated, $\si$ is presumed to be written as
 $$\si=(\si_1, \ldots, \si_{\ell(\si)}) \hspace{3mm} \textrm{ with }\hspace{3mm} \si_1 \geq \cdots \geq \si_{\ell(\si)}.$$ 
  To make an emphasis, if $\tau_k$ is another partition, it is simply $(\tau_{k1}, \ldots, \tau_{k\ell(\tau_k)}).$
  
 \item We say that $|\si|=n$ if $\si_1+\cdots +\si_{\ell(\si)}=n$.
 
 \item Let $\vect{\al}=(\al_1, \ldots, \al_{\ell(\si)})$ be an $\ell(\si)$-tuple of cohomology classes associated to $\si$ so that we may form a cohomology-weighted partition 
$$\si(\vect{\al}):=\si_1(\al_1)\cdots \si_{\ell(\si)}(\al_{\ell(\si)}).$$
The group $\Aut(\si(\vect{\al}))$ is defined to be the group of permutations on the set $\{1, \ldots, \ell(\si) \}$ fixing $\left( (\si_1, \al_1),\ldots, (\si_{\ell(\si)}, \al_{\ell(\si)})\right)$.
Moreover, let $\Aut(\si)$ be the group $\Aut(\si(\vect{\al}))$ when all entries of $\vect{\al}$ are identical.
 
 \item Let $\fz_{\si}=|\Aut(\si)|\prod_{i=1}^{\ell(\si)}\si_{i}$ be the order of the centralizer of any permutation of cycle type $\si$.
 
 \item We let $$(2)=(2, 1^{n-2}) \hspace{5mm} \textrm{and} \hspace{5mm} 1=(1^n)$$ be partitions of length $n-1$ (for $n \geq 2$) and length $n$ respectively.

 \end{enumerate}
 
\item Suppose $\la_i$ is a partition of a nonnegative integer $n_i$ for $i=1, \ldots, s$, and $n=\sum_{i=1}^s n_i$. The $s$-tuple $(\la_1, \ldots, \la_s)$ of partitions is denoted by $$\wt{\la}.$$ We also use the same symbol for the class \eqref{ptcls} in the localized equivariant Chen--Ruan cohomology of $[\Sym^n(S)]$.

\end{enumerate}

\section{Symmetric Product Stack}\label{s: Sym}
\subsection{Some definitions}\label{s: DefSym}
As mentioned earlier, $S$ always indicates a smooth toric surface. Let $n$ be a positive integer. For any nonempty subset $N$ of $\{1, \ldots, n \}$, let $\fS_N$ be the symmetric group on $N$ and 
$$S^N=\{(s_i)_{i\in N}: \spa s_i\textrm{'s are elements of } S \},$$ 
a set of $|N|$-tuples of elements of $S$.
For simplicity of notation, we denote by $\fS_n$ the group $\fS_{\{1,\ldots, n\}}$ and by $S^n$ the set $S^{\{1,\ldots, n\}}$.

The $n$-fold symmetric product $\Sym^n(S)=S^n/\fS_n$ of $S$ is the coarse moduli scheme of the quotient stack $[\Sym^n(S)]$ defined as follows: 
\begin{enumerate}
\item[(a)] An object over $U$ is a pair $(p \co P\to U, \, f \co P \to S^n)$, where $p$ is a principal $\fS_n$-bundle, and $f$ is an $\fS_n$-equivariant morphism.

\item[(b)] Suppose that $(p^\prime \co P^\prime \to U^\prime, f^\prime  \co P^\prime \to S^n)$ is another object. A morphism from $(p^\prime, f^\prime)$ to $(p, f)$ is a Cartesian diagram
$$
\begin{CD}
P^\prime @> \al >> P\\
@V p^\prime VV @VV p V\\
U^\prime @> \be >> U\\
\end{CD}
$$
such that $f^\prime = f \circ \al $. 
\end{enumerate}
The symmetric product stack $[\Sym^n(S)]$ is an orbifold with the natural atlas $S^n \to [S^n/\fS_n]$.

\subsection{Chen--Ruan cohomology}\label{s: ChenRuan}
Let $\I[\Sym^n(S)]$ be the stack of cyclotomic gerbes in $[\Sym^n(S)]$ (see Abramovich--Graber--Vistoli \cite{AGV}). It is isomorphic to a disjoint union of orbifolds 
$$\coprod_{[g]\in C} [S^n_g/\overline{C(g)}],$$ 
where $C$ is the set of conjugacy classes of $\fS_n$, $C(g)$ is the centralizer of $g$, $\overline{C(g)}$ is the quotient group $C(g)/\langle g \rangle$, and $S^n_g$ is the $g$-fixed locus of $S^n$.
Obviously, the connected components of $\I[\Sym^n(S)]$ can be labeled with the partitions of $n$. If $[g]$ is the conjugacy class corresponding to the partition $\la$, we may write
$$ 
S(\la)=S^n_g/C(g) \hspace{5mm} \textrm{and} \hspace{5mm}  \overline{S(\la)}=S^n_g/\overline{C(g)}.
$$
The component $[S^n/\fS_n]$ is called the untwisted sector while all other components of the stack $\I[\Sym^n(S)]$ are called twisted sectors. 

The Chen--Ruan cohomology
$$\Aorb([\Sym^n(S)])$$ 
is defined to be the cohomology $\Ast(\I[\Sym^n(S)])$ of the stack $\I[\Sym^n(S)]$.
Thus, it is simply
$$\bigoplus_{[g]\in C}\Ast(S^n_g/C(g))=\bigoplus_{[g]\in C}\Ast(S^n_g)^{C(g)}.$$
(For any orbifold $\cX$ with coarse moduli space $X$, we identify $\Ast(\cX)$ with $\Ast(X)$ via the pushforward $c_*\co \Ast(\cX) \to \Ast(X)$ defined by $c_*([\cV])=\frac{1}{s}[c(\cV)]$, where $\cV$ is a closed integral substack of $\cX$, and $s$ is the order of the stabilizer of a generic geometric point of $\cV$.)

Additionally, for $\al\in A^i(S(\la))$, the orbifold (Chow) degree of $\al$ is defined to be $i+\age(\la)$, where $\age(\la)=n-\ell(\la)$ is the age of $S(\la)$. In other words,
$$A_{\mathrm{orb}}^*([\Sym^n(S)])=\bigoplus_{|\la|=n} A^{*-\age(\la)}(S(\la)).$$

As $S$ carries a $\bT$-action, we may put the above cohomologies into an equivariant context by considering $\bT$-equivariant cohomologies as the quotient schemes and orbifolds discussed above clearly inherit $\bT$-actions from $S$.

\subsection{Bases and fixed-point classes} 
In this section, we exhibit a basis for the equivariant Chen--Ruan cohomology of the symmetric product stack $[\Sym^n(S)]$. It will be helpful for the determination of extended 3-point functions and for setting up our desired SYM-HILB correspondence later.

Given a partition $\la$ of $n$, we would like to construct a basis for $\AT(S^n_{g})^{C(g)}$, where $g\in \fS_n$ has cycle type $\la$.

The permutation $g$ has a cycle decomposition (i.e., a product of disjoint cycles including 1-cycles)
$$
g=g_1 \cdots g_{\ell(\la)}
$$
with $g_i$ being a $\la_i$-cycle. 
For every $i$, let $N_i$ be the smallest subset of $\{1,\ldots, n\}$ such that $g_i \in \fS_{N_i}$. Thus, $|N_i|=\la_i$, and $\coprod_{i=1}^{\ell(\la)}  N_i=\{1,\ldots, n\}$. It is clear that
$$
S^n_g=\prod_{i=1}^{\ell(\la)} S^{N_i}_{g_{i}} \hspace{5mm} \textrm{and} \hspace{5mm} \spa S^{N_i}_{g_i}\cong S.
$$

To the partition $\la$, we associate an $\ell(\la)$-tuple $\vect{\eta}=(\eta_1, \ldots, \eta_{\ell(\la)})$ of classes in $\AT(S)$. Let us put
\begin{equation}\label{basis}
g(\vect{\eta}) =\frac{1}{|\Aut(\la(\vect{\eta}))|\prod_{i=1}^{\ell(\la)} \la_i} \sum_{h\in C(g)}\bigotimes_{i=1}^{\ell(\la)} g_i^h(\eta_i) \in \AT(S^n_{g})^{C(g)}.
\end{equation}
This expression requires some explanation:
\begin{enumerate}
\item[(a)] $g_i^h$ denotes $h^{-1}g_i h$, and the class $g_i^h(\eta_i)$ is the pullback of the class $\eta_i$ on $S$ by the obvious isomorphism $S^{h^{-1} N_i h}_{h^{-1}g_i h} \xrightarrow{\cong} S$. (Note that we multiply permutations from left to right.) 

\item[(b)] Two classes 
$\bigotimes_{i=1}^{\ell(\la)} g_i^{h_1}(\eta_i)$ and $\bigotimes_{i=1}^{\ell(\la)} g_i^{h_2}(\eta_i)$
on the space $S^n_g$ coincide for some $h_1, h_2 \in C(g)$, and a straightforward verification shows that each $\bigotimes_{i=1}^{\ell(\la)} g_i^h(\eta_i)$ repeats as many as 
$$|\Aut(\la(\vect{\eta}))|\prod_{i=1}^{\ell(\la)} \la_i$$ 
times. Hence, $(|\Aut(\la(\vect{\eta}))|\prod_{i=1}^{\ell(\la)} \la_i)^{-1}$ is a normalization factor to ensure that no repetition occurs in \eqref{basis}.

\item[(c)] If $g=k_1 \cdots k_{\ell(\la)}$ is another cycle decomposition with $k_i$ being a $\la_i$-cycle, then there exists $h\in C(g)$ such that
$$\bigotimes_{i=1}^{\ell(\la)} k_i(\eta_i)=\bigotimes_{i=1}^{\ell(\la)} g_i^h(\eta_i).$$
Thus, the expression \eqref{basis} is independent of the cycle decomposition of $g$.
\end{enumerate}

Given a basis $\fB$ for $\AT(S)$. The classes $g(\vect{\eta})$'s, with $\eta_i$'s running over all elements of $\fB$, form a basis for $\AT(S^n_g)^{C(g)}$. (Note that if $\hat{g}$ is another permutation of cycle type $\la$, the classes $g(\vect{\eta})$ and $\hat{g}(\vect{\eta})$ are identical in the equivariant Chen--Ruan cohomology $\ATorb([\Sym^n(S)])$.) 

From now on, we use the notation 
$$
\la(\vect{\eta})
$$
for the class $g(\vect{\eta})$ in \eqref{basis}.
The classes $\la(\vect{\eta})$'s, ranging over all partitions $\la$ of $n$ and all $\eta_i \in \fB$, give a basis for $\ATorb([\Sym^n(S)])$.

We are going to work with the $\bT$-fixed point basis $\{[x_1], \ldots, [x_s] \}$.
Let $\la_i$ be a partition of a nonnegative integer $n_i$ for $i=1, \ldots, s$, and let $n=\sum_{i=1}^s n_i$. We denote the class
\begin{equation}\label{ptcls}
\la_{11}([x_1])\cdots \la_{1\ell(\la_1)}([x_1]) \cdots \la_{s1}([x_s])\cdots \la_{s\ell(\la_s)}([x_s])
\end{equation}
by
$$\wt{\la}:=(\la_1, \ldots, \la_s).$$
This class corresponds to a $\bT$-fixed point, which we denote by 
$$[\wt{\la}],$$
in the sector indexed by the partition 
$(\la_{11}, \ldots, \la_{1 \ell(\la_1)}, \ldots, \la_{s1}, \ldots, \la_{s \ell(\la_s)})$. So we refer to $\wt{\la}$'s as $\bT$-fixed point classes, which form a basis for the localized equivariant Chen--Ruan cohomology $\ATorb([\Sym^n(S)])_\fm$.

In the case of the affine plane $U_i$, the point $x_i$ is the unique $\bT$-fixed point. The $\bT$-fixed point class $$\la_{i1}([x_i])\cdots \la_{i\ell(\la_i)}([x_i]) \in \ATorb([\Sym^{n_i}(U_i)])_\fm$$
is denoted by
$$\wt{\la_i}.$$

Moreover, for every $\bT$-fixed point class $\wt{\la}$,
\begin{equation} \label{tangent weight}
\eT(T_{[\wt{\la}]}\bar{I}[\Sym^n(S)])=\prod_{k=1}^s {\eT(T_{x_k} U_k)}^{\ell(\la_k)}
\end{equation}
where $T_{[\wt{\la}]}\bar{I}[\Sym^n(S)]$ is the tangent space to $[\Sym^n(S)]$ at the $\bT$-fixed point $[\wt{\la}]$, $T_{x_k} U_k$ is the tangent space to $U_k$ at $x_k$, and $\eT( \, \cdot \, )$ indicates the $\bT$-equivariant Euler class.

\subsection{Extended three-point functions} \label{s: 3-extftn}
\subsubsection{Twisted stable maps and extended Gromov--Witten invariants} \label{s: OrbGW}
For any positive integers $m$ and $n$, we denote by
$$\overline{M}_{0,m}([\Sym^n(S)])$$ 
the moduli space parametrizing genus zero, $m$-pointed, twisted stable map 
$$f \co (\cC,\cP_1,\ldots, \cP_m)\to [\Sym^n(S)]$$ 
of degree zero. Note, in particular, that $\cP_i\cong \mathcal{B} \mu_{r_i}$ is the classifying stack of the cyclic group $\mu_{r_i}$ of $r_i$-th roots of unity for some positive integer $r_i$. The map $f$ is representable and comes equipped with the ordinary stable map
$$f_c \co (C, P_1,\ldots, P_m) \to \Sym^n(S),$$ 
where $C:=c(\cC)$ is the coarse curve of $\cC$, and $P_i=c(\cP_i)$ for $i=1, \ldots, m$. For more information on twisted stable maps, see Chen--Ruan \cite{CR2} or Abramovich--Graber--Vistoli \cite{AGV}.

Also, for every $i=1, \ldots, m$, there is an evaluation map
$$
\ev_i \co \M_{0,m}([\Sym^n(S)]) \to \I[\Sym^n(S)]
$$
which takes $[f \co (\cC,\cP_1,\ldots, \cP_m) \to [\Sym^n(S)]]$ to $[f \big\vert_{\cP_i} \co \cP_i \to [\Sym^n(S)]]$.
For any partitions $\tau_1,...,\tau_m$ of $n$ and nonnegative integer $d$, let 
$$
\M([\Sym^n(S)],\tau_1,...,\tau_m; d)=\bigcap_{i=1}^m \ev_i^{-1}([\overline{S(\tau_i)}])\cap \bigcap_{j=1}^d \ev_{m+j}^{-1}([\overline{S((2))}])
$$ 
be an open and closed substack of the moduli space $\M_{m+d}([\Sym^n(S)])$. (If $n=1$ and $d \geq 1$, the substack is taken to be empty.)

Given any $\al_i \in \ATorb([\Sym^n(S)])_\fm$ for $i=1, \ldots, m$, the $m$-point, degree zero, extended Gromov--Witten invariant 
$$
\langle \al_1, \ldots, \al_m \rangle_{d}^{[\Sym^n(S)]}
$$
is defined by
\begin{equation}\label{definv}
\frac{1}{d!} \sum_{|\tau_1|, \ldots, |\tau_m|=n}\int_{[\M([\Sym^n(S)],\tau_1,\ldots,\tau_m; \, d)]^{\vir}_{\bT}}\ev_1^*(\al_1)\cdots \ev_m^*(\al_m)
\end{equation}
(cf. Bryan--Graber \cite{BrG}). Here $[\hspace{3mm}]^{\vir}_{\bT}$ represents the equivariant virtual class. Note that the moduli spaces $\M([\Sym^n(S)],\tau_1,\ldots,\tau_m; d)$ are not necessarily compact, but the above definition makes sense when $\al_i$'s are $\bT$-fixed point classes because the loci of twisted stable maps meeting $\bT$-fixed points are compact. If $\al_i$'s are any classes in $\ATorb([\Sym^n(S)])_\fm$, \eqref{definv} can be defined by writing each $\al_i$ in terms of $\bT$-fixed point classes and by linearity. Another interpretation is to treat \eqref{definv} as a sum of residue integrals over $\bT$-fixed connected components of the spaces $\M([\Sym^n(S)],\tau_1,\ldots,\tau_m; d)$ via the virtual localization formula of Graber and Pandharipande \cite{GP}. In both intepretations, \eqref{definv} lies in the field $\bQ(t_1,t_2)$. 

We are primarily interested in 3-point extended invariants, and we encode them in a generating function: For $\al_1$, $\al_2$, $\al_3 \in \ATorb([\Sym^n(S)])_\fm$, set
$$
\langle \al_1, \al_2, \al_3 \rangle^{[\Sym^n(S)]}(u)=\sum_{d=0}^\infty \langle \al_1, \al_2, \al_3 \rangle_{d}^{[\Sym^n(S)]}u^d.
$$
We refer to these functions as extended 3-point functions. 

Note that the Chen--Ruan product \cite{CR1} is given by the (non-extended) 3-point orbifold Gromov--Witten invariants in degree zero (i.e., by setting $u=0$), and so the extended 3-point functions should provide more enumerative information than the Chen--Ruan product.

\subsubsection{Fixed locus}
We would like to describe the $\bT$-fixed locus of the space $\M([\Sym^n(S)], \tau_1, \ldots, \tau_m; d)$.
Suppose $$[f \co (\cC,\cP_1,\ldots, \cP_m, \cQ_1, \ldots, \cQ_d)\to [\Sym^n(S)]]$$ is an arbitrary $\bT$-fixed point of $\M([\Sym^n(S)],\tau_1,...,\tau_m; d)$. It naturally comes with the following diagram
\begin{equation}\label{diagram}
\begin{CD} 
P_{\cC} @> f^\prime >> S^n\\
@V  VV @VV \pi V\\
\cC @>f >> [\Sym^n(S)]\\
@V c VV @VV c V\\
C @>f_c >> \Sym^n(S).
\end{CD}
\end{equation}
Here $\pi$ is the natural atlas, and $P_{\cC}$ is the fiber product $\cC\times_{[\Sym^n(S)]}S^n$. Note that $P_{\cC}$ is a scheme by the representability of $f$. In addition, taking $f^\prime$ modulo $\fS_{n-1}$ and composing with the $n$-th projection, we have a map 
$$\wtd{f} \co \wtd{C} \to S^{\bT} \subseteq S$$
where $\wtd{C}=P_{\cC}/\fS_{n-1}$ is in fact a degree $n$ admissible cover of $C$, which is ramified over the markings $c(\cP_1),\ldots, c(\cP_m)$ with ramification profiles $\tau_1, \ldots, \tau_m$ and is simply ramified over the last $d$ markings $c(\cQ_1),\ldots, c(\cQ_d)$.
For ease of explanation, the points $c(\cQ_1),\ldots, c(\cQ_d)$ (resp. $\cQ_1, \ldots, \cQ_d$) are referred to as simple marked points of $C$ (resp. $\cC$).

In what follows, we use the notation 
$$\wt{a}$$ 
for the $s$-tuple $(a_1, \ldots, a_s)$ of nonnegative integers. Moreover, we define
$$
|\wt{a}|=\sum_{i=1}^s a_i.
$$

Suppose that $\wt{\si_1}, \ldots, \wt{\si_m}$ are any $s$-tuples of partitions such that the partition $\tau_j$ admits a decomposition $\tau_j=(\si_{j1}, \ldots, \si_{js})$ for $j=1, \ldots, m$ and $|\si_{1k}|=\cdots=|\si_{mk}|$ for $k=1, \ldots, s$, and that $\wt{a}$ is any $s$-tuple of nonnegative integers with $|\wt{a}|=d$. Let
$$
\fM(\wt{\si_1}, \ldots, \wt{\si_m}; \wt{a})
$$ 
be the locus in the space $\M([\Sym^n(S)], \tau_1, \ldots, \tau_m; d)$ which parametrizes $\bT$-fixed, $(m+d)$-pointed, twisted stable maps $f \co \cC \to [\Sym^n(S)]$ of degree zero satisfying the following properties: 
\begin{enumerate}
\item[(i)] The associated admissible cover $\wtd{C}$ may be written as a disjoint union 
$$\coprod_{k=1}^s \wtd{C}_k,$$
where each $\wtd{C}_k$, if nonempty, is contracted by $\wtd{f}$ to $x_k$. (Note that $\wtd{C}_k$ is possibly empty or disconnected. Empty sets are included just for simplicity of notation.)

\item[(ii)] For every $k=1, \ldots, s$, the cover $\wtd{C}_k \to C$ is ramified with monodromy
$$\si_{1k}, \ldots, \si_{mk}$$
above the first $m$ markings of $C$ and is ramified with monodromy 
$$(2)^{a_k}, 1^{d-a_k}$$ 
above the simple markings of $C$.
\qed
\end{enumerate}

From the discussion following Diagram \eqref{diagram}, we see that the (disjoint) union of the spaces $\fM(\wt{\si_1}, \ldots, \wt{\si_m}; \wt{a})$'s is $\M([\Sym^n(S)], \tau_1, \ldots, \tau_m; d)^\bT$.

Suppose $m \geq 3$. Let $n_k=|\si_{ik}|$ for $i=1, \ldots, m$. Then we have a natural morphism 
$$\phi \co \fM(\wt{\si_1}, \ldots, \wt{\si_m}; \wt{a}) \to \prod_{k=1}^s \M([\Sym^{n_k} U_k], \si_{1k}, \ldots, \si_{mk}; a_k)^\bT$$
defined as follows:
Let $[f \co \cC \to [\Sym^n(S)]]$ be an element of $\fM(\wt{\si_1}, \ldots, \wt{\si_m}; \wt{a})$. 
For each $k=1, \ldots, s$, we stabilize the target of the cover $\wtd{C}_k \to C$ and the domain accordingly (by forgetting those simple markings of $C$ which are {\em not} branched points of the cover $\wtd{C}_k \to C$). The output is the following setting
$$
\begin{CD} 
\wtd{C_{\st_k}} @>  >> \{x_i \}\\
@V   VV \\
C_{\st_k}\\
\end{CD}
$$
with the vertical map being an admissible cover of degree $n_k$ over the stabilization $C_{\st_k}$ of $C$. It gives rise to a $\bT$-fixed twisted stable map, which we denote by
\begin{equation}\label{stabilize}
f_k \co \cC_{\st_k} \to [\Sym^{n_k}(U_k)].
\end{equation}
The map $f_k$ represents a $\bT$-fixed point of the space $\M([\Sym^{n_k} U_k], \si_{1k}, \ldots, \si_{mk}; a_k)$. 
We then take $\phi([f])=([f_1], \ldots, [f_s])$.

Now we let $m=3$, in which case the morphism $\phi$ is surjective, and both its source and target have dimension $d$.
For any $\bT$-fixed connected component $F(\wt{a})$ of $\fM(\wt{\si_1}, \wt{\si_2}, \wt{\si_3}; \wt{a})$, we let
$$
F_k(\wt{a})={\pi_k \circ \phi (F(\wt{a}))}
$$
where $\pi_k$ is the $k$-th projection.
The collection $\prod_{k=1}^s F_k(\wt{a})$'s form a complete set of $\bT$-fixed connected components of
$\prod_{k=1}^s  \M([\Sym^{n_k} U_k], \si_{1k},  \si_{2k}, \si_{3k}; a_k)$.

Before proceeding, we fix some notation.
For every $\bT$-fixed connected component $F$ of the moduli space $\M([\Sym^n(S)],\tau_1, \tau_2, \tau_3; d)$, let
$$\iota_F \co F \to \M([\Sym^n(S)],\tau_1, \tau_2, \tau_3; d)$$
be the natural inclusion, and let
$$N_F^{\vir}$$ 
be the virtual normal bundle to $F$. 

\subsubsection{The product formula}
We want to investigate the invariant $\langle \wt{\la}, \wt{\mu}, \wt{\nu} \rangle^{[\Sym^n(S)]}_d$ for any $\bT$-fixed point classes $\wt{\la}, \wt{\mu}, \wt{\nu}$ and any nonnegative integer $d$. It is clearly zero if the condition
\begin{equation}\label{condition}
n_k:=|\la_k|=|\mu_k|=|\nu_k| \textrm{ for each } k=1,\ldots,s
\end{equation}
fails to hold. For the rest of this section, we assume that the $\bT$-fixed point classes $\wt{\la}, \wt{\mu}, \wt{\nu}$ satisfy Condition \eqref{condition}, and that $\sum_{i=1}^s n_i=n$.

\begin{prop}\label{prod-stack}
For any nonnegative integer $d$,
\begin{equation} \label{prod-formula}
\langle \wt{\la}, \wt{\mu}, \wt{\nu} \rangle^{[\Sym^n(S)]}_d=\sum_{a_1+\cdots+a_s=d} \, \prod_{k=1}^s \langle \wt{\la_k}, \wt{\mu_k}, \wt{\nu_k} \rangle^{[\Sym^{n_k}(U_k)]}_{a_k}.
\end{equation}
\end{prop}
\proof
The only $\bT$-fixed connected components which may contribute to the 3-point extended invariant
\begin{equation}\label{three point}
\langle \wt{\la}, \wt{\mu}, \wt{\nu} \rangle^{[\Sym^n(S)]}_d
\end{equation}
are the components of $\fM(\wt{\la}, \wt{\mu}, \wt{\nu}; \wt{a})$'s, where $|\wt{a}|=d$. Precisely, \eqref{three point} is given by
$$
\frac{1}{d!} \sum_{|\wt{a}|=d} \sum_{F(\wt{a})}\int_{F(\wt{a})}\frac{\iota^*_{F(\wt{a})}(\ev_1^*(\wt{\la})\cdot \ev_2^*(\wt{\mu})\cdot \ev_3^*(\wt{\nu}))}{\eT(N^\vir_{F(\wt{a})})},
$$
where $F(\wt{a})$ runs over all components of $\fM(\wt{\la}, \wt{\mu}, \wt{\nu}; \wt{a})$. 

For any component $F(\wt{a})$ of $\fM(\wt{\la}, \wt{\mu}, \wt{\nu}; \wt{a})$ and any point $[f] \in F(\wt{a})$, we have
\begin{eqnarray*}
\eT(H^i(\cC, f^*T[\Sym^n(S)]))=\phi^* \bigotimes_{k=1}^s \eT(H^i(\cC_{\st_k}, f_k^*T[\Sym^{n_k} U_k]))
\end{eqnarray*}
for $i=0, 1$. (Here we follow the notation of \eqref{stabilize}.) As a result,
$$
\eT(N^\vir_{F(\wt{a})})=\phi^* \bigotimes_{k=1}^s \eT(N^\vir_{F_k(\wt{a})}).
$$
Moreover,
$$
\eT(T_{[\wt{\si}]}\bar{I}[\Sym^n(S)])=\prod_{k=1}^s \eT(T_{[\wt{\si_k}]}\bar{I}[\Sym^{n_k}(U_k)])
$$
for each $\bT$-fixed point class $\wt{\si}$ (see \eqref{tangent weight}). 
This forces
$$\iota^*_{F(\wt{a})}(\ev_1^*(\wt{\la})\cdot \ev_2^*(\wt{\mu})\cdot \ev_3^*(\wt{\nu}))=\prod_{k=1}^s \iota^*_{F_k(\wt{a})}(\ev_1^*(\wt{\la_k})\cdot \ev_2^*(\wt{\mu_k}) \cdot \ev_3^*(\wt{\nu_k})).$$
Hence the contribution of $\fM(\wt{\la}, \wt{\mu}, \wt{\nu}; \wt{a})$ to \eqref{three point} equals
$$
\frac{1}{a_1!\cdots a_s!} \sum_{F(\wt{a})} \prod_{k=1}^s \int_{F_k(\wt{a})}\frac{\iota^*_{F_k(\wt{a})}(\ev_1^*(\wt{\la_k})\cdot \ev_2^*(\wt{\mu_k})\cdot \ev_3^*(\wt{\nu_k}))}{\eT(N^\vir_{F_k(\wt{a})})},
$$
where the prefactor accounts for the distribution of simple marked points. The sum is nothing but
\begin{eqnarray} \label{sum}
\prod_{k=1}^s \frac{1}{a_k!}  \sum_{F_k(\wt{a})}  \int_{F_k(\wt{a})}\frac{\iota^*_{F_k(\wt{a})}(\ev_1^*(\wt{\la_k})\cdot \ev_2^*(\wt{\mu_k})\cdot \ev_3^*(\wt{\nu_k}))}{\eT(N^\vir_{F_k(\wt{a})})}.
\end{eqnarray}
Since $F_k(\wt{a})$'s run through the $\bT$-fixed connected components of $\M([\Sym^{n_k} U_k], \wt{\la_k}, \wt{\mu_k}, \wt{\nu_k}; a_k)$, we obtain the equality \eqref{prod-formula} by summing \eqref{sum} over all nonnegative integers $a_1, \ldots, a_s$ that add up to $d$.
\qed

Put another way, we have the following product formula.
\begin{prop}\label{break-s} 
The extended 3-point function $\langle \wt{\la}, \wt{\mu}, \wt{\nu} \rangle^{[\Sym^n(S)]}(u)$ is given by
$$
\prod_{k=1}^s \langle \wt{\la_k}, \wt{\mu_k}, \wt{\nu_k} \rangle^{[\Sym^{n_k}(U_k)]}(u).
$$
Moreover, every extended 3-point function is a rational function in $t_1, t_2, e^{iu}$, where $i^2=-1$.
\end{prop}
\proof
The first statement is immediate from Proposition \ref{prod-stack}. Also, the 3-point function $\langle \wt{\la}, \wt{\mu}, \wt{\nu} \rangle^{[\Sym^n(S)]}(u)$ lies in $\bQ(t_1, t_2, e^{iu})$ because each $\langle \wt{\la_k}, \wt{\mu_k}, \wt{\nu_k} \rangle^{[\Sym^{n_k}(U_k)]}(u)$ is an element of $\bQ(t_1, t_2, e^{iu})$ (see Okounkov--Pandharipande \cite{OP} and Bryan--Graber \cite{BrG}). As every equivariant Chen--Ruan cohomology class is a $\bQ(t_1, t_2)$-linear combination of $\bT$-fixed point classes, the second statement follows.
\qed

From the above results, we find that the 3-point extended invariants of $[\Sym^n(S)]$ in degree zero are completely determined by the invariants of the symmetric product stacks of $\bC^2$. 

The orbifold Poincar\'e pairing $\langle \cdot \, | \, \cdot \rangle$ on $\ATorb([\Sym^{n_k}(U_k)])_{\fm}$ is determined by
\begin{equation}\label{e:orbifoldpairing}
\langle \wt{\la_k}\, | \,\wt{\mu_k} \rangle=\de_{\la_k, \mu_k}\frac{(L_k R_k)^{\ell(\la_k)}}{\fz_{\la_k}},
\end{equation}
where $\de_{\la_k, \mu_k}$ stands for the Kronecker delta.
The argument of Proposition \ref{prod-stack} may be applied to determine the orbifold Poincar\'e pairing on $\ATorb([\Sym^n(S)])_{\fm}$. 
\begin{prop}\label{p:orbifoldpair}
Let $\wt{\si}$ and $\wt{\tau}$ be any $\bT$-fixed point classes in $\ATorb([\Sym^n(S)])_{\fm}$. Then the pairing $\langle \wt{\si}\, | \,\wt{\tau} \rangle$ is
\begin{eqnarray*}
\prod_{k=1}^s  \frac{(L_k R_k)^{\ell(\si_k)}}{\fz_{\si_k}}
\end{eqnarray*}
if $\wt{\si}=\wt{\tau}$ and zero otherwise. Thus, the $\bT$-fixed point classes are orthogonal with respect to the orbifold Poincar\'e pairing. \qed
\end{prop}

\section{Hilbert scheme of points}\label{s: Hilb}
\subsection{Fixed-point basis and Nakajima basis}
\paragraph{Fixed-point basis.}
As shown by Ellingsrud and Str\o mme \cite{ES}, there is a one-to-one correspondence between partitions of $n$ and $\bT$-fixed points of $\Hilb^n(\bC^2)$. Since $U_i\cong \bC^2$, for every partition $\la$, the corresponding $\bT$-fixed point $$\la(x_i)\in \Hilb^{|\la|}(U_i)$$ can be described as follows: Suppose $\bC[u,v]$ is the coordinate ring of $U_i$. Then $\la(x_i)$ is the zero-dimensional closed subscheme of $U_i$ whose ideal $\cI_{\la(x_i)}$ is
$$(u^{\la_1},v u^{\la_2},\ldots, v^{\ell(\la)-1}u^{\la_{\ell(\la)}},v^{\ell(\la)}).$$ 
The point $\la(x_i)$ is supported at $x_i$ and is mapped to $|\la| \, [x_i]\in \Sym^{|\la|}(U_i)$ by the Hilbert--Chow morphism.

The action of $\bT$ on $S$ lifts to $\Hilb^n(S)$, and the fixed locus $\Hilb^n(S)^{\bT}$ is isolated. Each $\bT$-fixed point of $\Hilb^n(S)$ is supported in $\{x_1,\ldots,x_s\}$ and may be expressed as a sum
$$
\la_1(x_1)+\cdots+\la_s(x_s)
$$
for some partitions $\la_i$'s satisfying $\sum_{i=1}^s |\la_i|=n$ (by ``the sum'' we mean the disjoint union of $\la_i(x_i)$'s). For $\wt{\la}=(\la_1,\ldots,\la_s)$, write
$$
I_{\wt{\la}}=[\la_1(x_1)+\cdots+\la_s(x_s)],
$$
which are $\bT$-fixed point classes of $\Hilb^n(S)$ and form a basis for the localized equivariant cohomology $\AT(\Hilb^n(S))_{\fm}$. 

\paragraph{Nakajima basis.}
Another important basis for $\AT(\Hilb^n(S))$ is the Nakajima basis, which we now describe. For further details, see Grojnowski \cite{Gro}, Nakajima \cite{N1,N2}, Vasserot \cite{V}, or Li--Qin--Wang \cite{LQW2}.

Given a partition $\la$ of $n$ and an associated $\ell(\la)$-tuple $\vect{\eta}=(\eta_1,\ldots, \eta_{\ell(\la)})$ of classes in $\AT(S)$. We define
$$\fa_{\la}(\vect{\eta})=\frac{1}{|\Aut(\la(\vect{\eta}))|}\prod_{i=1}^{\ell(\la)}\frac{1}{\la_i}\spa \fp_{-\la_i}(\eta_i)|0\rangle,$$
where $|0\rangle=1 \in A_{\bT}^0(S^{[0]})$, and $\fp_{-\la_i}(\eta_i)\co \AT(S^{[k]})\to A_{\bT}^{*+\la_i-1+\mathrm{deg}(\eta_i)/2}(S^{[k+\la_i]})$ are Heisenberg creation operators. Note that we also denote the class $\fa_{\la}(\vect{\eta})$ by 
$$\fa_{\la_1}(\eta_1)\cdots \fa_{\la_{\ell(\la)}}(\eta_{\ell(\la)}).$$ 

Choose a basis $\fB$ for $\AT(S)$. The classes
$\fa_{\la}(\vect{\eta})$'s, running over all partitions $\la$ of $n$ and all $\eta_i \in \fB$,
give a basis for $\AT(\Hilb^n(S))$. They are referred to as the Nakajima basis with respect to $\fB$. 

Moreover, we may work with the Nakajima basis with respect to the $\bT$-fixed point basis $\{[x_1],\ldots,[x_s]\}$.
Let $\la_1, \ldots, \la_s$ be partitions of $n_1, \ldots, n_s$ respectively, and let $n=\sum_{i=1}^s n_i$. We define
$$\fa_{\wt{\la}}=\fa_{\la_{11}}([x_1])\cdots \fa_{\la_{1\ell(\la_1)}}([x_1]) \cdots \fa_{\la_{s1}}([x_s])\cdots \fa_{\la_{s\ell(\la_s)}}([x_s]).$$
Set $\ell(\wt{\la})=\sum_{i=1}^s \ell(\la_i)$. The Chow degree of $\fa_{\wt{\la}}$ is
\begin{equation}\label{Chow deg}
\sum_{i=1}^s (n_i-\ell(\lambda_i))+ 2\ell(\wt{\la})=(n-\ell(\wt{\la}))+2\ell(\wt{\la})=n+\ell(\wt{\la}).
\end{equation}
In the case of $U_i$, the class
$$\fa_{\la_{i1}}([x_i])\cdots \fa_{\la_{i\ell(\la_i)}}([x_i]) \in \AT(\Hilb^{n_i}(U_i))$$
is denoted by 
$$\fa_{\wt{\la_i}}.$$
The Poincar\'e pairing $\langle \cdot \, | \, \cdot \rangle$ on $\AT(\Hilb^{n_i}(U_i))_{\fm}$ is determined by the formula
\begin{equation}\label{e:hilbertpairing}
\langle \fa_{\wt{\la_i}} | \fa_{\wt{\mu_i}} \rangle=\delta_{\la_i,\mu_i} \frac{(-1)^{|\la_i|-\ell(\la_i)}(L_i R_i)^{\ell(\la_i)}}{\fz_{\la_i}},
\end{equation} 
where $|\la_i|=|\mu_i|=n_i$.

\subsection{Comparison to symmetric functions}\label{s:comparison}
Let $p_i(z)=\sum_{k=1}^\infty z_k^i$ be the $i$-th power sum for each positive integer $i$, and let $p_0(z)=1$. For partitions $\la$, write
$$
p_{\la}(z)=\frac{1}{|\Aut(\la)|}\prod_{i=1}^{\ell(\la)} \frac{1}{\la_i}p_{\la_i}(z),
$$
which form a basis for the algebra $\mathrm{\La}_{\, \bQ(t_1,t_2)}$ of symmetric functions in the variables $\{ z_k \}_{k=1}^\infty$ over $\bQ(t_1,t_2)$. Let $\al_i=-R_i/L_i$ for $i=1, \ldots, s$. For any partition $\mu$, we denote by
$$J_{\mu}^{(\al_i)}(z)$$
the integral form of the Jack symmetric function indexed by $\mu$ and $\al_i$. Note that for each $\al_i$, the symmetric functions $J_{\mu}^{(\al_i)}(z)$'s provide an orthogonal basis for $\mathrm{\La}_{\, \bQ(t_1,t_2)}$; see Macdonald \cite{Mac}.

Furthermore, we want to understand the relationship between the $\bT$-fixed point basis and the Nakajima basis with respect to $\{[x_1], \ldots, [x_s]\}$. This amounts to relating the Jack symmetric functions and the power sums. Indeed, the $\bT$-fixed point class $I_{(\mu_1,\ldots,\mu_s)}$ is identified with 
$$\bigotimes_{i=1}^s L_i^{|\mu_i|} J_{\mu_i}^{(\al_i)}(z(i)),$$
while the Nakajima basis element $\fa_{(\la_1,\ldots,\la_s)}$ is identified with 
$$\bigotimes_{i=1}^s L_i^{\ell(\la_i)}p_{\la_i}(z(i)).$$ 
For more details, see \cite{LQW2,N1,V}.

For $i=1,\ldots, s$, let $\la_i$ be a partition of $n_i$. As the fixed-point classes $[\mu_i(x_i)]$'s (with $|\mu_i|=n_i$) span $\AT(\Hilb^{n_i}(U_i))_{\fm}$, we can write 
$$\fa_{\wt{\la_i}}=\sum_{|\mu_i|=n_i} c_{\la_i,\mu_i}[\mu_i(x_i)]$$ 
for some $c_{\la_i,\mu_i}\in \bQ(t_1, t_2)$. By the above identifications,
\begin{equation}\label{decomp}
\fa_{(\la_1,\ldots,\la_s)}=\sum c_{\la_1,\mu_1}\cdots c_{\la_s,\mu_s} \, I_{(\mu_1,\ldots,\mu_s)}
\end{equation}
where the sum is over all partitions $\mu_i$ of $n_i$, $i=1,\ldots,s$.

\subsection{Extremal three-point functions}\label{three point GW invariants}\label{s: 3pt}
\subsubsection{Extremal invariants and quantum corrected cohomology} \label{s: extrGW}
Let $\rho \co \Hilb^n(S) \to \Sym^n(S)$ be the Hilbert--Chow morphism as before. If $n \geq 2$, the kernel of the induced homomorphism
$$\rho_* \co A_1(\Hilb^n(S);\bZ)\to A_1(\Sym^n(S);\bZ)$$
is one-dimensional and is generated by an effective rational curve class $\be_n$ that is Poincar\'e dual to $-\fa_{2}(1)\fa_1(1)^{n-2}$. For any integers $k \geq 1$ and $d \geq 0$, we let $$\M_{0,k}(\Hilb^n(S),d)$$ 
be the moduli space parametrizing stable maps from genus zero, $k$-pointed, nodal curves to $\Hilb^n(S)$ of degree $d\be_n$. (If $n=1$ and $d>0$, we take $\M_{0,k}(\Hilb^n(S),d)$ to be empty.)

Let $e_i \co \M_{0,k}(\Hilb^n(S),d)\to \Hilb^n(S)$ be the evaluation map at the $i$-th marked point, and let $\al_i \in  \AT(\Hilb^n(S))_{\fm}$ be any equivariant class for $i=1, \ldots, k$.
Although $\Hilb^n(S)$ is not necessarily compact, the $k$-point extremal Gromov--Witten invariant 
$$
\langle \al_1, \ldots, \al_k \rangle^{\Hilb^n(S)}_{d} = \int_{[\M_{0,k}(\Hilb^n(S), \, d)]^{\mathrm{vir}}_{\bT}} e_1^*(\al_1)\cdots e_k^*(\al_k)
$$
is well-defined by a similar treatment to \eqref{definv}.

We will explore the following extremal 3-point functions of $\Hilb^n(S)$:
$$\langle \al_{1}, \al_{2}, \al _{3}\rangle^{\Hilb^n(S)}(q)=\sum_{d=0}^{\infty}\langle \al_{1},\al_2,\al _{3}\rangle_{d}^{\Hilb^n(S)} q^d.$$

The quantum corrected product $\qcr$ for $\Hilb^n(S)$ is defined by the above generating function with the specialization $q=-1$; see Ruan \cite{R}. Precisely, the product $a\qcr b$ of any two classes $a$, $b$ in $\AT(\Hilb^n(S))_{\fm}$ is defined to be the unique element satisfying
$$
\langle a \qcr b \, | \, c \rangle = \langle a,b,c \rangle^{\Hilb^n(S)}(q)\big\vert_{q=-1}, \spa \forall \, c \in  \AT(\Hilb^n(S))_{\fm}.
$$
Here $\langle \cdot \, | \, \cdot \rangle$ denotes the Poincar\'e pairing on $\AT(\Hilb^n(S))_{\fm}$. 
The equivariant cohomology $\AT(\Hilb^n(S))_{\fm}$ endowed with the multiplication $\qcr$ is referred to as the equivariant quantum corrected cohomology ring of $\Hilb^n(S)$. 

On the other hand, the cup product $\cup$ for $\Hilb^n(S)$ is given by the rule:
$$\langle a \cup b \, | \, c \rangle =\langle a,b,c \rangle^{\Hilb^n(S)}(q)\big\vert_{q=0}.$$
Note that the right-hand side is simply the degree zero invariant $\langle a,b,c \rangle^{\Hilb^n(S)}_0$.

\subsubsection{The product formula}\label{s: Hilb.4}
In this section, we would like to determine the extremal Gromov--Witten invariants of the Hilbert scheme $\Hilb^n(S)$ in Nakajima basis elements $\fa_{\wt{\la}}$, $\fa_{\wt{\mu}}$, $\fa_{\wt{\nu}}$. We will see that these invariants may be expressed in terms of the invariants of the Hilbert schemes of points on $\bC^2$. 
First of all, we have the following. 

\begin{prop} \label{p:vanishing}
For every nonnegative integer $d$, the invariant $\langle \fa_{\wt{\la}},\fa_{\wt{\mu}}, \fa_{\wt{\nu}}\rangle_{d}^{\Hilb^n(S)}$ does not vanish only if
\begin{equation} \label{e:mono}
|\la_i|=|\mu_i|=|\nu_i| \spa \textrm{ for each }i=1,\ldots,s. 
\end{equation}
\end{prop}
\proof
Let us suppress the superscript $\Hilb^n(S)$ for the moment.
By \eqref{decomp}, the invariant $\langle \fa_{\wt{\la}},\fa_{\wt{\mu}}, \fa_{\wt{\nu}}\rangle_{d}$ is a $\bQ(t_1, t_2)$-linear combination of the invariants $\langle I_{\wt{\si}}, I_{\wt{\tau}}, I_{\wt{\theta}} \rangle_{d}$ where
\begin{equation} \label{e:mono*}
|\si_i|=|\la_i|, \hspace{3mm} |\tau_i|=|\mu_i|, \hspace{3mm} |\theta_i|=|\nu_i|, \hspace{3mm} \forall i=1,\ldots,s. 
\end{equation}

Suppose that $\langle \fa_{\wt{\la}}, \fa_{\wt{\mu}}, \fa_{\wt{\nu}}\rangle_{d} \not=0$. Then by the preceding discussion,
\begin{equation} \label{e:fixed pt inv}
\langle I_{\wt{\si}}, I_{\wt{\tau}}, I_{\wt{\theta}} \rangle_{d}\not= 0
\end{equation}
for some $\wt{\si}, \wt{\tau}, \wt{\theta}$ satisfying \eqref{e:mono*}. By virtual localization, the invariant $\langle I_{\wt{\si}}, I_{\wt{\tau}}, I_{\wt{\theta}} \rangle_{d}$ is expressed as a sum of residue integrals over $\bT$-fixed connected components. Because of  \eqref{e:fixed pt inv}, there is a $\bT$-fixed component, say $\G$, which gives a nonvanishing contribution to the sum.

Now let $[f \co (C, p_{\wt{\si}}, p_{\wt{\tau}}, p_{\wt{\theta}}) \to \Hilb^n(S)] \in \G$ be an arbitrary point. As $f$ is of degree $d\be_n$, the image of the composite morphism $\rho \circ f$ must be a point of $\Sym^n(S)$. Moreover, for $\wt{\eta}=\wt{\si}, \wt{\tau}, \wt{\theta}$, the point $f(p_{\wt{\eta}})$ must be the $\bT$-fixed point $\eta_1(x_1)+\cdots + \eta_s(x_s)$ due to the evaluation condition imposed by the insertion $I_{\wt{\eta}}$. This means that
$$\rho \circ f (p_{\wt{\eta}})=\sum_{i=1}^s |\eta_i| \, [x_i].$$
Note that all $\rho \circ f (p_{\wt{\eta}})$'s are the same point. Hence, $|\si_i|=|\tau_i|=|\theta_i|$ for $i=1,\ldots, s$, and we conclude that \eqref{e:mono} holds.
\qed  

It remains to study the 3-point invariants 
$$\langle \fa_{\wt{\la}},\fa_{\wt{\mu}},\fa_{\wt{\nu}}\rangle_{d}^{\Hilb^n(S)}$$
for all $\wt{\la}$, $\wt{\mu}$, $\wt{\nu}$ satisfying Condition \eqref{e:mono}: $n_i:=|\la_i|=|\mu_i|=|\nu_i|$ for each $i=1, \ldots, s$, and $\sum_{i=1}^s n_i=n$. We fix such $n_i$ and partitions $\la_i$, $\mu_i$, $\nu_i$ of $n_i$ throughout the remainder of this section.

Let 
$$U=\Hilb^{n_1}(U_1)\times \cdots \times \Hilb^{n_s}(U_s) \hspace{5mm}\textrm{and } \hspace{5mm}  P=\rho^{-1}(n_1[x_1]+\cdots+n_s[x_s]).$$
In fact, 
$$P\cong \rho_1^{-1}(n_1[x_1])\times \cdots \times \rho_s^{-1}(n_s[x_s])\subseteq U,$$ 
where $\rho_i \co  \Hilb^{n_i}(S) \to \Sym^{n_i}(S)$ is the Hilbert--Chow morphism for every $i$. (Note that in case $n_i=0$, $\Hilb^{n_i}(U_i)$ and $\rho_i^{-1}(n_i[x_i])$ will be dropped from the above products.)
Suppose $N=\{i=1, \ldots, s \, | \, n_i\geq 1 \}$. As $\rho_i^{-1}(n_i[x_i])$ is irreducible and has complex dimension $n_i-1$ for $i \in N$, $P$ is then irreducible and has dimension $n-|N|$.

Let $\xi=\mu_1(x_1)+ \cdots +\mu_s(x_s) \in \Sn$ and $\xi_U=\mu_1(x_1)\times \cdots \times\mu_s(x_s) \in U$, where $\mu_i$ is a partition of $n_i$ for each $i$. These two points are identical in $P$. In addition, we have the identification $T_{\xi_U}U=T_\xi \Sn$ of tangent spaces. Indeed,
\begin{eqnarray*}
T_{\xi_U}U &=&\bigoplus_{i \in N}\textrm{Hom}_{\cO_{U_i}}(\cI_{\mu_i(x_i)},\cO_{\mu_i(x_i)})\\
&=&\bigoplus_{i \in N}\textrm{Hom}_{\cO_{S,x_i}}(\cI_{\xi, x_i},\cO_{\xi, x_i})\\
&=&\textrm{Hom}_{\cO_{S}}(\cI_{\xi},\cO_{\xi})\\
&=&T_{\xi}\Sn.
\end{eqnarray*}

Denote by $\io_P, \jm_P$ the inclusion of $P$ into $\Hilb^n(S)$ and $U$ respectively. We have a simple lemma.
\begin{lem}\label{l:jacksu}
$\io_P^*(\fa_{\wt{\la}})=\jm_P^*(\fa_{\wt{\la_1}} \otimes \cdots \otimes \fa_{\wt{\la_s}}).$
\end{lem}
\proof
By \eqref{decomp}, it suffices to show that
\begin{equation}\label{e:fix}
\io_P^*I_{(\mu_1,\ldots,\mu_s)}=\jm_P^*([\mu_1(x_1)]\otimes \cdots \otimes [\mu_s(x_s)]).
\end{equation}
Let $\xi$ and $\xi_U$ be the points as in the discussion preceding the lemma.
We see that the left-hand side of \eqref{e:fix} is given by
\begin{eqnarray*}
\sum_{\eta \in P^{\bT}} \frac{{i_\eta}_* (\io_P\circ i_\eta)^*I_{(\mu_1,\ldots,\mu_s)}}{\eT(T_\eta P)}
= \frac{\eT(T_{\xi} \Sn)}{\eT(T_\xi P)},
\end{eqnarray*}
where $i_{\eta} \co \{\eta\} \to P$ is the natural inclusion. Similarly, the right-hand side of \eqref{e:fix} coincides with
$$\frac{\eT(T_{\xi_U} U)}{\eT(T_{\xi_U} P)}.$$ 
Thus, the equality \eqref{e:fix} follows from $T_{\xi_U}U=T_\xi \Sn$.
\qed

\begin{prop}\label{p:breaking}
For every nonnegative integer $d$,
$$
\langle \fa_{\wt{\la}},\fa_{\wt{\mu}},\fa_{\wt{\nu}}\rangle^{\Hilb^n(S)}_{d}=\sum_{d_1+\cdots+d_s=d} \, \prod_{i=1}^s \langle \fa_{\wt{\la_i}},\fa_{\wt{\mu_i}},\fa_{\wt{\nu_i}} \rangle^{\Hilb^{n_i}(U_i)}_{d_i}.
$$
\end{prop}
\proof
To determine the 3-point invariant $\langle \fa_{\wt{\la}},\fa_{\wt{\mu}},\fa_{\wt{\nu}}\rangle_{d}^{\Hilb^n(S)}$, we only need to consider those connected components of $\N(\Hilb^n(S),d)^{\bT}$ whose images under the map $\rho\circ e_i$ are the point 
$$n_1[x_1]+\cdots+n_s[x_s]$$ 
for every $i=1,2,3$. Observe that any $\bT$-fixed stable map $f$ that represents a point of any of these components factors through $P$. 
Hence, we may calculate $\langle \fa_{\wt{\la}},\fa_{\wt{\mu}},\fa_{\wt{\nu}}\rangle^{\Hilb^n(S)}_{d}$ over the connected components of
$$\coprod_{d_1+\cdots+d_s=d}\N(P,(d_1,\ldots,d_s))^{\bT}.$$

We denote by $\g \co \G \to \N(\Hilb^n(S),d)$ the natural inclusion for every $\bT$-fixed connected component $\G$.
However, we write $\G_{d_1,\ldots,d_s}$ for $\G$ whenever $\G$ is contained in the space $\N(P,(d_1,\ldots,d_s))$. On the other hand, we also note that $\G_{d_1,\ldots,d_s}$'s form a complete set of $\bT$-fixed connected components of $\N(U,(d_1,\ldots,d_s))$. The following diagram summarizes their relationships:
$$
\xymatrix{&\G_{d_1,\ldots,d_s} \ar[ld]_{\g_U} \ar[d]^{\g_P} \ar[rd]^{\g} \\ \N(U,(d_1,\ldots,d_s))  &\N(P,(d_1,\ldots,d_s)) \ar@{_{(}->}[l] \ar@{^{(}->}[r] &\N(\Hilb^n(S),d)}
$$
where $d=d_1+\cdots+d_s$, and $\g_P$ and $\g_U$ are the natural inclusions.

We would like to show that
\begin{equation}\label{prod}
\langle \fa_{\wt{\la}},\fa_{\wt{\mu}},\fa_{\wt{\nu}}\rangle^{\Hilb^n(S)}_{d}=\sum_{d_1+\cdots+d_s=d}\langle \otimes_{i=1}^s \fa_{\wt{\la_i}}, \otimes_{i=1}^s \fa_{\wt{\mu_i}}, \otimes_{i=1}^s \fa_{\wt{\nu_i}} \rangle^{U}_{(d_1,\ldots,d_s)}.
\end{equation}
Denote by $\bae_i \co \N(P,(d_1,\ldots,d_s)) \to P$ the evaluation map at the $i$-th marked point. The invariant $\langle \fa_{\wt{\la}},\fa_{\wt{\mu}},\fa_{\wt{\nu}}\rangle^{\Hilb^n(S)}_{d}$ is
\begin{eqnarray*}
\sum_{d_1+\cdots+d_s=d} \, \sum_{\G_{d_1,\ldots,d_s}} \int_{\G_{d_1,\ldots,d_s}} \frac{\g_P^*( \bae_1^*\io_P^*(\fa_{\wt{\la}})\cdot \bae_2^*\io_P^*(\fa_{\wt{\mu}})\cdot \bae_3^*\io_P^*(\fa_{\wt{\nu}}))}{\eT(N_{\G_{d_1,\ldots,d_s}}^\vir)}
\end{eqnarray*}
where $N_{\G_{d_1,\ldots,d_s}}^{\vir}$ is the virtual normal bundle to $\G_{d_1,\ldots,d_s}$ in $\N(\Hilb^n(S),d)$.
On the other hand, we find that $\langle \otimes_{i=1}^s \fa_{\wt{\la_i}}, \otimes_{i=1}^s \fa_{\wt{\mu_i}}, \otimes_{i=1}^s \fa_{\wt{\nu_i}} \rangle^{U}_{(d_1,\ldots,d_s)}$ is 
$$\sum_{\G_{d_1,\ldots,d_s}} \int_{\G_{d_1,\ldots,d_s}}\frac{\g_P^*(\bae_1^*\jm_P^*(\otimes_{i=1}^s \fa_{\wt{\la_i}})\cdot \bae_2^*\jm_P^*(\otimes_{i=1}^s \fa_{\wt{\mu_i}})\cdot \bae_3^*\jm_P^*(\otimes_{i=1}^s \fa_{\wt{\nu_i}}))}{\eT(N_{\G_{d_1,\ldots,d_s}, U}^\vir)}$$
where $N_{{\G_{d_1,\ldots,d_s}}, U}^\vir$ is the virtual normal bundle to $\G_{d_1,\ldots,d_s}$ in $\N(U,(d_1,\ldots,d_s))$.
By Lemma \ref{l:jacksu}, it is given by
$$
\sum_{\G_{d_1,\ldots,d_s}}\int_{\G_{d_1,\ldots,d_s}}\frac{\g_P^* ( \bae_1^*\io_P^*(\fa_{\wt{\la}})\cdot \bae_2^*\io_P^*(\fa_{\wt{\mu}})\cdot \bae_3^*\io_P^*(\fa_{\wt{\nu}}))}{\eT(N_{{\G_{d_1,\ldots,d_s}}, U}^\vir)}.
$$
Hence, \eqref{prod} is a consequence of the equality on the inverse equivariant Euler classes of virtual normal bundles:
\begin{equation}\label{e:eulernormal}
\frac{1}{\eT(N_{\G_{d_1,\ldots,d_s}}^\vir)}=\frac{1}{\eT(N_{{\G_{d_1,\ldots,d_s}}, U}^\vir)} \spa \textrm{  for each } \G_{d_1,\ldots,d_s}.
\end{equation}

Now we prove \eqref{e:eulernormal}. Suppose $[f \co \Si \to P]$ is any point of the $\bT$-fixed component $\G_{d_1,\ldots,d_s}$. For 
$$
X=S^{[n]} \hspace{5mm} \textrm{or} \hspace{5mm} \spa X=U,
$$
in order to verify \eqref{e:eulernormal}, we need only examine the infinitesimal deformations of $f$ with the source curves fixed. In fact, it suffices to check that
$$
\frac{\eT(H^0(\Si_v,f^*TX)^{\mov})}{\eT(H^1(\Si_v,f^*TX)^{\mov})} \hspace{5mm} \textrm{and} \hspace{5mm} \frac{\eT(H^0(\Si_e,f^*TX)^{\mov})}{\eT(H^1(\Si_e,f^*TX)^{\mov})}
$$ are independent of $X$ for every connected contracted component $\Si_v$ and noncontracted irreducible component $\Si_e$. Here, for example, $H^0(\Si_v,f^*TX)^{\mov}$ denotes the moving part of $H^0(\Si_v,f^*TX)$. 

The first independence holds due to the fact that $\Si_v$ is of genus $0$ and $T_{f(\Si_v)}S^{[n]}=T_{f(\Si_v)}U$. Thus, it remains to justify the second independence. Since $\Si_e\cong \bP^1$, $f^*TX$ is a direct sum of line bundles over $\Si_e$, i.e., $f^* TX=\bigoplus_{i=1}^{2n}\cO_{\Si_e}(\ell_i^X)$ for some integers $\ell_i^X$'s, and we also have
\begin{equation}\label{split}
\frac{\eT(H^0(\Si_e,f^*TX)^{\mov})}{\eT(H^1(\Si_e,f^*TX)^{\mov})}=\prod_{i=1}^{2n}\frac{\eT(H^0(\Si_e,\cO_{\Si_e}(\ell_i^X))^{\mov})}{\eT(H^1(\Si_e,\cO_{\Si_e}(\ell_i^X))^{\mov})}.
\end{equation}
Note that the $\bT$-action on $f(\Si_e)$ (independent of $X$) induces a $\bT$-action on $\Si_e$ and hence actions on $\cO_{\Si_e}(\ell_i^X)$'s. Let $p_0$ and $p_\infty$ be the points of $\Si_e$ which correspond respectively to the points $0$ and $\infty$ of $\bP^1$ via $\Si_e\cong \bP^1$.
Suppose that for $j=0, \infty$, $\bT$ acts on $\cO_{\Si_e}(1)\vert_{p_j}$ with weight $w_j$. Then the $\bT$-weight of $\cO_{\Si_e}(\ell_i^X)\vert_{p_j}$ is $\ell_i^X w_j$. This means that $T_{f(p_j)}X$ comes with $\bT$-weights $\ell_1^X w_j,\ldots,\ell_{2n}^X w_j$. As $T_{f(p_j)}S^{[n]}=T_{f(p_j)}U$, $(\ell_1^{S^{[n]}} w_j,\ldots, \ell_{2n}^{S^{[n]}} w_j)$ and $(\ell_1^U w_j,\ldots, \ell_{2n}^U w_j)$ are equal after a suitable reordering. By \eqref{split}, this shows the second independence and ends the proof of \eqref{e:eulernormal}.

On the other hand, the equality
$$
\langle \otimes_{i=1}^s \fa_{\wt{\la_i}}, \otimes_{i=1}^s \fa_{\wt{\mu_i}}, \otimes_{i=1}^s \fa_{\wt{\nu_i}} \rangle^{U}_{(d_1,\ldots,d_s)}=\prod_{i=1}^s \langle \fa_{\wt{\la_i}},\fa_{\wt{\mu_i}},\fa_{\wt{\nu_i}}\rangle^{\Hilb^{n_i}(U_i)}_{d_i}
$$
holds. This is clear for $(d_1,\ldots,d_s)=(0,\ldots,0)$; in general, the above equality follows from the product formula of Behrend \cite{Beh} in equivariant context. Combining it with \eqref{prod}, we obtain the proposition.
\qed

\begin{prop} \label{p:product formula}
The 3-point function $\langle \fa_{\wt{\la}},\fa_{\wt{\mu}},\fa_{\wt{\nu}}\rangle^{\Hilb^n(S)}(q)$ lies in $\bQ(t_1,t_2,q)$ and is given by
$$
\prod_{i=1}^s\langle \fa_{\wt{\la_i}},\fa_{\wt{\mu_i}},\fa_{\wt{\nu_i}}\rangle^{\Hilb^{n_i}(U_i)}(q).
$$
\end{prop}
\proof
Each $\langle \fa_{\wt{\la_i}},\fa_{\wt{\mu_i}},\fa_{\wt{\nu_i}}\rangle^{\Hilb^{n_i}(U_i)}(q)$ is a rational function in $t_1, t_2, q$ (cf. \cite{OP}).
By Proposition \ref{p:breaking},
\begin{eqnarray*}
\langle \fa_{\wt{\la}},\fa_{\wt{\mu}},\fa_{\wt{\nu}}\rangle^{\Hilb^n(S)}(q)&=&\sum_{d=0}^{\infty}(\sum_{d_1+\ldots+d_s=d}\prod_{i=1}^s \langle \fa_{\wt{\la_i}},\fa_{\wt{\mu_i}},\fa_{\wt{\nu_i}}\rangle^{\Hilb^{n_i}(U_i)}_{d_i})q^d\\
&=&\prod_{i=1}^s\sum_{d_i=0}^{\infty}\langle \fa_{\wt{\la_i}},\fa_{\wt{\mu_i}},\fa_{\wt{\nu_i}}\rangle^{\Hilb^{n_i}(U_i)}_{d_i}q^{d_i}\\
&=&\prod_{i=1}^s\langle \fa_{\wt{\la_i}},\fa_{\wt{\mu_i}},\fa_{\wt{\nu_i}} \rangle^{\Hilb^{n_i}(U_i)}(q),
\end{eqnarray*}
which is a rational function in $t_1, t_2, q$ as well.
\qed

From Eq. \eqref{e:hilbertpairing}, we obtain the following result about the Poincar\'e pairing on the localized equivariant cohomology $\AT(\Hilb^n(S))_{\fm}$.
\begin{prop} \label{p:hilbertpair}
Suppose $\wt{\si}$ and $\wt{\tau}$ are any $s$-tuples of partitions such that both $\fa_{\wt{\si}}$ and $\fa_{\wt{\tau}}$ are classes in $\AT(\Hilb^n(S))_{\fm}$. Then the pairing $\langle \fa_{\wt{\si}}| \fa_{\wt{\tau}}\rangle$ is
$$
\prod_{i=1}^s \frac{(-1)^{|\si_i|-\ell(\si_i)}(L_i R_i)^{\ell(\si_i)}}{\fz_{\si_i}}
$$
if $\wt{\si}=\wt{\tau}$ and zero otherwise.
In other words, the Nakajima basis $\{ \fa_{\wt{\si}} \}$ is orthogonal with respect to the Poincar\'e pairing.
\qed
\end{prop}

\section{SYM-HILB correspondence}\label{s:3}
\subsection{Proof of the main theorem}
In this section, we give a SYM-HILB correspondence that relates the theories of the symmetric product stack $[\Sym^n(S)]$ and the Hilbert scheme $\Hilb^n(S)$ for every positive integer $n$. 

Let $q=-e^{iu}$, where $i^2=-1$, and $K=\bQ(i, t_1,t_2)$. Define
$$L \co \ATorb([\Sym^n(S)])\otimes_{\bQ[t_1,t_2]} K \to \AT(\Hilb^n(S))\otimes_{\bQ[t_1,t_2]} K$$ 
by
$$
L(\wt{\la})=(-i)^{\age(\wt{\la})}\fa_{\wt{\la}}.
$$
Here $\age(\wt{\la})=\sum_{i=1}^s \age(\la_i)$. As we have a bijection between the bases on both sides, $L$ extends to a $K((u))$-linear isomorphism. However, extended 3-point functions of $[\Sym^n(S)]$ and extremal 3-point functions of $\Hilb^n(S)$ are elements of $K(q)$. So we may view $L$ as a $K(q)$-linear isomorphism. 
Note also that $\wt{\la}$ has orbifold Chow degree 
$$2\ell(\wt{\la})+\age(\wt{\la})=n+\ell(\wt{\la}),$$
which matches the Chow degree of $\fa_{\wt{\la}}$; see \eqref{Chow deg}.

We now verify Theorem \ref{main}, which is restated as follows.
\begin{Thm}\label{main theorem}
Under the substitution $q=-e^{iu}$, $L$ equates the extended 3-point functions of $[\Sym^n(S)]$ to the extremal 3-point functions of $\Hilb^n(S)$. More precisely, for any Chen--Ruan cohomology classes $\al_1, \al_2, \al_3$,
$$
\langle \al_1, \al_2, \al_3 \rangle^{[\Sym^n(S)]}(u)=\langle L(\al_1), L(\al_2), L(\al_3) \rangle^{\Hilb^n(S)}(q).
$$
Moreover, $L$ is an isometric isomorphism. 
\end{Thm}
\proof
For every $\bT$-fixed point class $\wt{\la}$, Proposition \ref{p:orbifoldpair} and Proposition \ref{p:hilbertpair} say that 
$$\langle \wt{\la}\, | \,\wt{\la} \rangle=(-1)^{\age(\wt{\la})}\langle \fa_{\wt{\la}}\, | \,\fa_{\wt{\la}} \rangle.$$
Hence, $\langle \wt{\la}\, | \,\wt{\la} \rangle=\langle L(\wt{\la})\, | \, L(\wt{\la}) \rangle$, and so $L$ is an isometric isomorphism.

The proof of the first assertion relies on the case of $\bC^2$. Indeed, Okounkov--Pandharipande and Bryan--Graber determine the structures of equivariant quantum cohomology rings of $\Hilb^n(\bC^2)$ and $[\Sym^n(\bC^2)]$ respectively. Also, these rings are related by the correspondence 
$$
L_{\bC^2} \co \ATorb([\Sym^n(\bC^2)])\otimes_{\bQ[t_1,t_2]} K(q) \to \AT(\Hilb^n(\bC^2))\otimes_{\bQ[t_1,t_2]} K(q)
$$
which sends $\mu_1([0])\cdots \mu_{\ell(\mu)}([0])$ to $(-i)^{\age(\mu)}\fa_{\mu_1}([0])\cdots \fa_{\mu_{\ell(\mu)}}([0])$. It follows from \cite[Theorem 3.11]{BrG} that  $L_{\bC^2}$ preserves 3-point functions. That is, for any Chen--Ruan cohomology classes $\de_1,\de_2,\de_3$ on the orbifold $[\Sym^n(\bC^2)]$, we have
\begin{equation}\label{BG}
\langle\de_1,\de_2,\de_3 \rangle^{[\Sym^n(\bC^2)]}(u)=\langle L_{\bC^2}(\de_1),L_{\bC^2}(\de_2),L_{\bC^2}(\de_3)\rangle^{\Hilb^n(\bC^2)}(q).
\end{equation}
As $U_k\cong \bC^2$, we denote the map corresponding to $L_{\bC^2}$ by 
$$L_{U_k} \co \ATorb([\Sym^n(U_k)])\otimes_{\bQ[t_1,t_2]} K(q) \to \AT(\Hilb^n(U_k))\otimes_{\bQ[t_1,t_2]} K(q).$$

To show the first assertion, it is enough to establish that for all $\wt{\la}, \wt{\mu}, \wt{\nu}$ satisfying Condition \ref{condition},
\begin{equation}\label{in fixed point}
\langle \wt{\la}, \wt{\mu}, \wt{\nu} \rangle^{[\Sym^n(S)]}(u)=\langle L({\wt{\la}}),L({\wt{\mu}}),L({\wt{\nu}})\rangle^{\Hilb^n(S)}(q).
\end{equation}
In fact, by Proposition \ref{break-s} and Eq. \eqref{BG}, the left-hand side of \eqref{in fixed point} equals
$$
\prod_{k=1}^s \langle L_{U_k}(\wt{\la_k}),L_{U_k}(\wt{\mu_k}),L_{U_k}(\wt{\nu_k})\rangle^{\Hilb^{n_k}(U_k)}(q),
$$
which is, by Proposition \ref{p:product formula}, given by
\begin{equation}\label{products}
\left(\prod_{k=1}^s (-i)^{\age(\la_k)+\age(\mu_k)+\age(\nu_k)} \right) \langle \fa_{\wt{\la}}, \fa_{\wt{\mu}}, \fa_{\wt{\nu}} \rangle^{\Hilb^n(S)}(q).
\end{equation}
As $\age(\wt{\si})=\sum_{k=1}^s \age(\si_k)$ for any fixed-point class $\wt{\si}$, \eqref{products} is the right-hand side of \eqref{in fixed point}. This completes the proof.
\qed

As explained earlier, this theorem also indicates that $L$ provides a ring isomorphism between the equivariant Chen--Ruan cohomology of $[\Sym^n(S)]$ and the equivariant quantum corrected cohomology of $\Hilb^n(S)$ after setting $u=0$ and $q=-1$.

\subsection{The cup product structure for the Hilbert scheme of points}
The upshot of Theorem \ref{main theorem} is that the 3-point degree zero invariants of $\Hilb^n(S)$ are expressible in terms of the invariants of $[\Sym^n(S)]$.  
\begin{cor} \label{precup}
Given cohomology-weighted partitions $\la_1(\vect{\eta_1})$, $\la_2(\vect{\eta_2})$, $\la_3(\vect{\eta_3})$ of $n$,
the degree zero invariant 
$$
\langle \fa_{\la_1(\vect{\eta_1})},\fa_{\la_2(\vect{\eta_2})},\fa_{\la_3(\vect{\eta_3})}\rangle^{\Hilb^n(S)}_0$$ is given by
$$
i^{\sum_{k=1}^3\age(\la_k)} \lim_{u \to +i \, \infty} \langle \la_1(\vect{\eta_1}), \la_2(\vect{\eta_2}), \la_3(\vect{\eta_3}) \rangle^{[\Sym^n(S)]}(u).
$$
\end{cor}
\proof
The corollary follows immediately from Theorem \ref{main theorem} as $q$ approaches $0$ whenever $u$ approaches $+i \, \infty$.
\qed

This is the precise statement of Corollary \ref{cup}. Since the ordinary cup product for $\Hilb^n(S)$ is defined by 3-point degree zero invariants, Corollary \ref{precup} says that the 3-point extended invariants of $[\Sym^n(S)]$ in degree zero completely determine the cup product for $\Hilb^n(S)$. This is not covered by CCRC and may shed new light on explicit calculations of the ordinary cohomology rings of the Hilbert schemes of points on $S$.

Finally, we remark that the relative Gromov--Witten theory of the threefold $S\times \bP^1$ is very close to the orbifold theory; see \cite{BrP,Ch,M}. Indeed, studying the relative 3-point invariants of $S\times \bP^1$ in class 
$(0,n) \in A_1(S \times \bP^1; \bZ)=A_1(S; \bZ) \oplus \bZ$
is tantamount to studying the extended 3-point functions of $[\Sym^n(S)]$. In other words, the relative theory may yield an alternative way to compute the cup product for $\Hilb^n(S)$. 

\paragraph{Acknowledgements.}
Many thanks are due to Tom Graber, who introduced me to the fascinating world of algebraic geometry. This paper has benefited greatly from his inspiring ideas and helpful suggestions. I am also pleased to thank the referee for useful comments. A portion of this paper was written under support from the National Science Council, Taiwan.

\phantomsection

\medskip
\noindent
{\it E-mail address}: keng@mail.ncku.edu.tw \\
Department of Mathematics,
National Cheng Kung University\\
Tainan 701,
Taiwan

\end{document}